\documentclass[sn-mathphys-num]{sn-jnl}

\usepackage{graphicx}%
\usepackage{multirow}%
\usepackage{amsmath,amssymb,amsfonts}%
\usepackage{amsthm}%
\usepackage{mathrsfs}%
\usepackage[title]{appendix}%
\usepackage{xcolor}%
\usepackage{textcomp}%
\usepackage{manyfoot}%
\usepackage{booktabs}%
\usepackage{algorithm}%
\usepackage{algorithmicx}%
\usepackage{algpseudocode}%
\usepackage{listings}%

\theoremstyle{thmstyleone}%

\theoremstyle{thmstyletwo}%

\theoremstyle{thmstylethree}%

\newcommand{\be}{\begin{equation}} 
\newcommand{\ee}{\end{equation}}
 
\newcommand{\R}{\mathbb{R}} 
\newcommand{\app}{\approx} 
\newcommand{\N}{\mathbb{N}} 
\newcommand{\A}{\mathcal{A}}

\newcommand{\TT}{\mathcal{T}} 
\newcommand{\CC}{\mathcal{C}}

\newcommand{\bg}{\begin} 
\newcommand{\en}{\end} 
 
\newcommand{\Om}{\Omega} 
\newcommand{\om}{\omega} 
\newcommand{\ep}{\epsilon} 
\newcommand{\de}{\delta} 
\newcommand{\si}{\sigma} 
 
\newcommand{\tx}{\text} 
\newcommand{\tit}{\emph} 
\newcommand{\es}{\emptyset} 
\newcommand{\ity}{\infty} 
\newcommand{\se}{\subseteq} 
\newcommand{\mo}{{-1}}

\newcommand{\Ra}{\Rightarrow}

\begin{document} 

\title[Cournot's principle for measure-theoretic probability]{Cournot's principle for measure-theoretic probability} 


\author{\fnm{Galvan} \sur{Bruno} }

\affil{\orgaddress{\city{Trento}, \postcode{38122}, \country{Italy}}}


\abstract{

The problem of relating the mathematics of probability theory to the empirical world of experiments has been debated for centuries. One of the oldest solutions proposed for this problem is a principle that states that an event with probability close to 1 nearly certainly occurs in a single trial of an experiment. This principle is now called \textit{Cournot' principle}.

Cournot's principle was first formulated in the context of classical probability, in which the probability of any event is given, and the \textit{product rule}, i.e., the rule that the probability that two events occur in two separated trials is the product of their probabilities, can be deduced. On the contrary, in the modern measure-theoretic approach to probability, probability measures and experiments are separate entities that must be related in an appropriate way, and the product rule cannot be deduced.

In this paper, a version of Cournot's principle suitable for measure-theoretic probability is proposed. Therefore, the principle is reformulated as a criterion for relating probability measures and experiments, and the product rule is explicitly stated.

In spite of the vagueness of the notions involved, the new version is formulated in a rigorous manner and an exact result, namely, that at most one probability measure can be related to an experiment, is rigorously proven.
}




\maketitle 

\section{Introduction} 

The problem of relating the mathematics of probability theory to the empirical world of experiments has been debated for centuries. One of the oldest solutions proposed for this problem is a principle that states that an event with probability close to 1 (which a mathematical property) nearly certainly occurs in a single trial of an experiment (which is an empirical property). This principle is now called \tit{Cournot' principle}. Cournot's principle, and in particular a version of it suited to measure-theoretic probability, is the subject of this paper.

Let us shortly recall the long history of Cournot's principle\footnote{The name Cournot's principle appears for the fist time after the Second World War \cite{shafer1}, so that most authors who mention it do not use this name.}. This principle was first formulated by Jacob Bernoulli in his \tit{Ars Conjectandi} (1713), to explain the correspondence between the probability and the relative frequency of an event. However, Augustin Cournot seems to have been the first to say explicitly that the whole empirical predictive power of classical probability derives from this principle (1843). Paul L{\'e}vy expressed the thesis that Cournot's principle is probability's only bridge to reality. \'Emile Borel called it ``the only law of chance'' (1940). In 1933 Kolmogorov proposed a version of the principle more suitable for measure-theoretic probability. See \cite{shafer1, shafer2007cournot} for an extensive presentation of Cournot's principle and of its history until the first half of the twentieth century. A notion closely related to Cournot's principle is that of \tit{typicality}\footnote{We recall that an event of a probability space is said to be \tit{typical} if its probability is close to $1$.}. In recent years, the role of typicality and Cournot's principle in statistical, quantum, and Bohmian mechanics has been studied by various authors, for example: \cite{galvan:4, seldonboltz, goldstein2012typicality, allori2020statistical, maudlin2020grammar, durr2021typicality,  lazarovici2023typicality}.

\vspace{3mm}
Let us now look in more detail at how Cournot's principle works. The principle has been initially formulated in the context of classical probability. We recall that the classical probability of an event is defined as the ratio between the number of outcomes composing the event and the number of possible outcomes. This definition is based on the \tit{principle of indifference}, according to which all outcomes are equally probable. Let us show now how the principle, in the context of classical probability, explains that the relative frequency of an event in a sufficiently long sequence of trials of an experiment is close to its probability.

Let $E$ denote a random experiment, and let $E^n$ denote the experiment whose trials are composed of $n$ trials of $E$. A trivial calculation shows that, if $P$ is the classical probability measure of $E$, then the $n$-fold product of $P$, denoted by $P^n$, is the classical probability measure of $E^n$. In fact, let us suppose that $E$ has $k$ possible outcomes, so that the probability $P(\{\om\})$ of any outcome $\om$ is $1/k$. Then $E^n$ has $k^n$ possible outcomes, so that the probability of any outcome $(\om_1, \ldots, \om_n)$ of $E^n$ is $1/k^n=P(\{\om_1\}) \cdots P(\{\om_n\})=P^n(\{\om_1 \times  \cdots \times \om_n\})$. Hereafter, this property is referred to as \tit{the product rule (for independent trials)}.

Now, let $\ep$ be a ``small'' positive number, $A$ be an event of $E$, and $S_n(A, \ep)$ be the event of $E^n$ composed of the sequences of outcomes of $E$ for which the relative frequency of $A$ differs from $P(A)$ less than $\ep$. Bernoulli's theorem (i.e., the weak law or large numbers) states that $P^n(S_n(A, \ep)) \to 1$ for $n \to \ity$, so that $P^n(S_n(A, \ep)) \app 1$ for $n$ large enough. Therefore, for such $n$, Cournot's principle applied to $E^n$ states that nearly certainly the event $S_n(A, \ep)$ occurs, i.e., nearly certainly the relative frequency of $A$ in the actual sequence is close to $P(A)$.

\vspace{3mm}
The first decades of the 20th century saw the decline of classical probability and the advent of measure-theoretic probability, which culminated in 1933 with the publication by Kolmogorov of his celebrated book on the foundations of probability theory \cite{kolmogorov1956foundations}. In measure-theoretic probability, the classical definition of probability is no longer available, and probability spaces and experiments are separate entities that must be related in an appropriate way. Therefore, the classical formulation of Cournot's principles is not suitable for measure-theoretic probability for essentially the following two reasons: (i) while the classical formulation predicts an empirical consequence from a given probability, namely, classical probability, in measure-theoretic probability the principle must be formulated as a criterion for establishing the correct association between a probability measure and an experiment; (ii) in the absence of the definition of classical probability, in measure-theoretic probability the product rule is no longer deducible, and, therefore, it must be explicitly stated.

In its book, Kolmogorov proposes a formulation of the principle more suitable for measure-theoretic probability. Roughly speaking, this formulation is the following \cite[\S 2]{kolmogorov1956foundations}: Given a random experiment $E$, under certain conditions (not specified by Kolmogorov) there is a probability measure $P$ such that: (a) nearly certainly the relative frequency of any event $A$ in a long sequence of trials is close to $P(A)$, and (b) if $P(A)\app 1$ for some event $A$, then nearly certainly $A$ occurs in a single trial of the experiment.

Let us shortly analyze Kolmogorov's formulation. It is reasonable to argue that Kolmogorov implicitly assumes that the probability measure $P$ satisfying points (a) and (b) is the correct probability measure to associate with the experiment $E$, so that Kolmogorov formulation can actually be considered as a criterion for relating probability measures and experiments. We note that point (b) holds only for the pair $(P, E)$ but not for a generic pain $(P^n, E^n)$ (see \cite[p. 91]{shafer1} for a discussion about this particular choice of Kolmogorov). This means that the product rule is not stated, and, therefore, Kolmogorov is forced to introduce point (a) to relate probability and relative frequency. In fact, we recall that in the previous reasoning relating probability and relative frequency, Cournot's principle is applied to the pair $(P^n, E^n)$, with $n$ large, rather than to $(P, E)$.

In this paper, another formulation of Cournot's principle suitable for measure-theoretic probability is proposed. The main changes of the new formulation compared to the classical one are: (1) the explicit definition of the empirical property of practical certainty, (2) the formulation of the principle as a criterion for relating probability measures and experiments (analogously to Kolmogorov's formulation), and (3) the explicit formulation of the product rule (which is essentially the requirement that Kolmogorov's point (b) holds for a generic experiment pair $(P^n, E^n$). Moreover, particular attention has been paid to formulating the new version in a form that is as mathematically rigorous as possible.

Considering probability measures and experiments as separate entities allows for the possible existence of three particular situations, namely: (i) that two probability measures are associated with the same experiment, (ii) that two different experiments are associated with the same probability measure, and (iii) that there are experiments that cannot be associated with any probability measure. In this paper, it is rigorously proven that the first situation does not occur.

The paper is organized as follows: In Section \ref{math}, some preliminary notions are introduced. In Section \ref{prace}, the empirical notion of practical certainty is defined. In Section \ref{proba}, the new formulation of the principle is presented. In Section \ref{part}, the three particular situations are presented. Section \ref{conc} concludes the paper.

\section{Preliminary notions \label{math}}

\subsection{Mathematical formalism}

The notions of sample space and event space (a $\si$-algebra of subsets of the sample space) are well known. Let $\Om$ and $\A$ denote a sample space and its event space, respectively. Then the symbol $\A^n$ denotes the event space of $\Om^n$ generated by the measurable rectangles of $\Om^n$.

Probability measures and probability spaces are also well-known notions. The symbols $P, P_1$, and $P_2$ always denote probability measures on the event space $\A$. The symbol $P^n$ denotes the $n$-fold product of $P$.

For $\de \in [0,1]$, let us define:
\be
\TT(P^n, \de):=\{A \in \A^n: P^n(A) \geq \de\}.
\ee

The following theorem is used in Section \ref{part} to prove that an experiment is governed by at most one probability measure:

\bg{theo} \label{th} For all $n \in \N$, let $\CC_n \se \A^n$ be a class that does not contain disjoint events. If 
\[
\TT(P^n_1, \de_1) \cup \TT(P^n_2, \de_2) \se \CC_n
\]
for some $\de_1, \de_2 < 1$ and all $n \in \N$, then $P_1 = P_2.$
\end{theo}
\bg{proof} Let us prove that $P_1 \neq P_2$ contradicts the hypotheses. Suppose that $P_1(A) \neq P_2(A)$ for some $A \in \A$. For $n \in \N$, let us define the function $f_n:\Om^n \to \R$ as follows:
\[
f_n(\om_1, \ldots, \om_n):= \frac{1}{n} \sum_{i=1}^n \chi_A(\om_i),
\]
where $\chi_A$ is the characteristic function of the event $A$. If $I$ is an interval of $\R$, the set $S_n[I]:=f_n^\mo[I]$ contains all the sequences $(\om_1, \ldots, \om_n)$ in which the relative frequency of $A$ belongs to $I$. This set belongs to $\A^n$ because $f_n$ is measurable, and, moreover, $S_n[I_1] \cap S_n[I_2]=\es$ if $I_1 \cap I_2 = \es$. Let $I_1$ and $I_2$ be two disjoint open intervals of $\R$ containing $P_1(A)$ and $P_2(A)$, respectively. From Bernoulli's theorem (the weak law of large numbers) it follows that 
\[
\lim_{n \to \ity} P_i^n(S_n[I_i]) =1, \tx{ where } i=1,2.
\]
See, for example, \cite{papoulis1990probability}. This implies that $S_m[I_i] \in \TT(P^m_i, \de_i)$ for $i=1, 2$ and some $m$ large enough. As a consequence, both disjoint events $S_m[I_1]$ and $S_m[I_2]$ belong to $\CC_m$, which contradicts the hypothesis that $\CC_m$ does not contain disjoint events.
\en{proof}

A simple corollary of this theorem is the following:

\bg{cor} If $\TT(P^n_1, \de_1)=\TT(P^n_2, \de_2)$ for some $\de_1, \de_2 \in (\tfrac{1}{2}, 1)$ and all $n \in \N$, then $P_1=P_2$.
\end{cor}
\bg{proof} The hypothesis implies that
\[
\TT(P^n_1, \de_1)  \cup \TT(P^n_2, \de_2) \se \TT(P^n_1, \de_1) \tx{ for all } n \in \N,
\]
and the classes $\{\TT(P^n_1, \de_1)\}_{n \in \N}$ satisfy the property required by Theorem \ref{th}.
\end{proof}
In other words, for any $\de \in (\tfrac{1}{2}, 1)$, the classes $\{\TT(P^n, \de)\}_{n \in \N}$ uniquely determine $P$.

\subsection{Experiments}

The empirical notions of random experiment and trial are well known. Every experiment is associated with an event space. The symbols $E, E_1$, and $E_2$ always denote experiments with event space $\A$. The symbol $E^n$ denotes the experiment with event space $\A^n$ consisting of $n$ repetitions of $E$, i.e., the experiment whose trials are composed of $n$ trials of $E$.

\subsection{Vague sets}

In the next section, an empirical definition of a practically certain event is proposed. This definition is unavoidably vague, and the corresponding set of practically certain events is vague as well.

To manage mathematically vague sets, we introduce the notion of instance: An \tit{instance} of a vague set is an exact set compatible with the definition of the vague set. Obviously, a vague set admits many instances, and the set of the instances is itself vague.

Hereafter, when the symbol denoting a vague set appears in a mathematical expression, we say that the vague set satisfies the expression, meaning that each of its instances satisfies the expression. Let us give an example.

Let $(\Om, \A, P)$ be a probability space. The class of \tit{typical} events is the vague set
\[
\TT:=\{A \in \A: P(A) \app 1\}, 
\]
where the expression $P(A) \app 1$ means $1-\ep \leq P(A) \leq 1$, with $0 \leq \ep \ll 1$.

From the implication 
\[
P(A) \app 1 \tx{ and } P(B) \geq P(A) \Ra P(B) \app 1,
\]
one easily deduces that the instances of $\TT$ are of the form
\[
\{A \in \A: P(A) \geq \de\} \tx{ or } \{A \in \A: P(A) > \de\},
\]
where $\de \leq 1$. The usual meaning we give to the expression $\ep \ll 1$ is certainly compatible with the two implications
\[
\ep \ll 1 \Ra \ep \leq \tfrac{2}{3}, \tx{ and } \ep \leq 10^{-30} \Ra \ep \ll 1,
\]
from which one deduces that (every instance of) $\TT$ satisfies the inclusions
\[
\TT(P, 1- 10^{-30}) \se \TT \se \TT(P, \tfrac{2}{3}).
\]
We note that, for example, the set $\TT(P, 1 - 10^{-40})$ is not an instance of $\TT$, because it excludes events that are certainly typical. 

Finally, from the above results, it is trivial to deduce that (every instance of) $\TT$ satisfies the following properties:
\bg{itemize}
\item[(a)] $\Om \in \TT$;
\item[(b)] if $A \in \TT$, $ B \in \A$, and $A \se B$, then $B \in \TT$;
\item[(c)] the set $\TT$ does not contain disjoint events.
\end{itemize}

\section{Practically certain events\label{prace}} 

Let us call \tit{practically certain} an event that occurs nearly certainly in a single trial of an experiment and nearly always in a sequence of trials\footnote{This terminology is borrowed from Kolmogorov \cite[\S 2]{kolmogorov1956foundations}. Bernoulli called these events \tit{morally certain}. The term \tit{almost certain} is generally used in everyday language; we prefer not to use this term here because, in mathematical probability, it has a different meaning, and this could generate confusion.}. Cournot's principle somehow evokes the property of practical certainty without explicitly declaring it. Therefore, the idea is to explicitly recognize practical certainty as an objective empirical property of some events of an experiment. Let us propose the following operational definition for this property:

\bg{defn}[Practical certainty] \label{pcd}  \hspace{1mm}
\bg{itemize}
\item[(a)] Practical certainty is defined operationally as follows: We select an event and then perform a long sequence of trials; the event is practically certain if and only if its relative frequency in the sequence is very close to $1$.
\item[(b)] Like any measurement, the above procedure may sometimes produce the wrong result. However, this does not prevent us from considering practical certainty as an objective property of some events.
\item[(c)] A natural consequence of the above definition is that if we select a practically certain event and perform a single trial of the experiment, then the event occurs nearly certainly in the trial.
\end{itemize}
\end{defn}

Let $\CC(E^n) \se \A^n$ denote the class of the practically certain events of the experiments $E^n$. Obviously $\CC(E^n)$ is a vague set of events. For all $n \in \N$, the class $\CC(E^n)$ satisfies the same properties as typical events:
\bg{itemize}
\item[(a)] $\Om^n \in \CC(E^n),$
\item[(b)] if $A \in \CC(E^n)$, $B \in \A^n$, and $A \se B$, then $B \in \CC(E^n)$,
\item[(c)] $\CC(E^n)$ does not contain disjoint events.
\end{itemize}

These properties cannot be rigorously deduced from the empirical Definition \ref{pcd}. However, they can be easily justified. Let us consider, for example, property (c): Let us choose two events $A$ and $B$ of $E^n$ and carry out a long series of trials of $E^n$. If $A$ is practically certain, its relative frequency is close to $1$; this implies that the relative frequency of $B$ is not close to $1$, and therefore, $B$ is not practically certain. Properties (a) and (b) can be justified even easier.

In conclusion, we assume/postulate that the classes $\{\CC(E^n)\}_{n \in N}$ satisfy the above properties.

We note that these properties do not exhaust the properties that a class of practically certain events should satisfy. For example, it is reasonable to require that the intersection of three practically certain events cannot be empty either. However, these properties, in particular property (c), are sufficient for our purposes.

\section{Cournot's principle for measure-theoretic probability\label{proba}} 

In this section, in successive steps, we present the new formulation of the principle.

Using the formalism previously introduced, the classical formulation of Cournot's principle reads:
\[
\TT(P, \de) \se \CC(E) \tx{ for } \de \app 1,
\]
where $P$ is the classical probability of the experiment $E$. A first attempt to adapt Cournot's principle to measure-theoretic probability is the following:
\[
\tx{A probability measure } P \tx{ \tit{governs} an experiment } E \tx{ if } \TT(P, \de) \se \CC(E) \tx{ for } \de \app 1.
\]
The verb ``governs'' may be replaced, for example, by ``is a probabilistic model of'', ``is associated with'', etc. We note that this formulation is similar to that of Kolmogorov (without point (a)). 

As previously said, in the absence of the classical definition of probability, the above formulation does not guarantee the product rule. In other words, it does not guarantee that if $P$ governs $E$, then $P^n$ governs $E^n$. The solution of Kolmogorov to this problem is to introduce point (a), while the solution proposed here corresponds to the requirement that point (b) holds for all pairs $(P^n, E^n)$.

In conclusion, the new formulation of Cournot's principle we propose is the following:

\bg{cp}\
A probability measure $P$ governs an experiment $E$ if there is $\de < 1$ such that
\be \label{cps}
\TT(P^n, \de) \se \CC(E^n) \tx{ for all } n \in \N.
\ee
\en{cp}

The requirement that conditions (\ref{cps}) hold for all $n \in \N$ is arguably an idealization, mainly because it requires that any number of trials is independent in a rigorous mathematical sense (see again \cite[p. 91]{shafer1}). Nevertheless, it is adopted in this paper because it allows us to highlight the logical structure of Cournot's principle and to deduce a rigorous result. Moreover, its study is certainly a prerequisite for a possible future study of more realistic situations.

We note that the condition ``for $\de \app 1$'' of the preliminary formulations has been replaced by the condition ``there is $\de < 1$'' in the final formulation. This change is motivated by the fact that it makes the formulation more precise. Moreover, even if the second condition seems more general, it actually implies the first one. In fact, suppose to the contrary that condition (\ref{cps}) is satisfied for some $\de \not \app 1$, and that there is an event $A \in \A$ such that $\de \leq P(A) \not \app 1$. The event $A$ is practically certain because it belongs to $\TT(P, \de)$ and, therefore, also to $\CC(E)$. On the other hand, from Bernoulli's theorem it follows that the long-run relative frequency of $A$ is close to $P(A) \not \app 1$, which contradicts the fact that $A$ is practically certain. In conclusion, condition (\ref{cps}) can be satisfied only if $\de \app 1$.

\section{Three particular situations \label{part}}

The new formulation of Cournot's principle highlights the possible existence of three particular situations.

The first situation occurs if two (or more) probability measures govern the same experiment. In formal terms, this happens if there are $\de_1, \de_2 < 1$ such that 
\bg{equation} \label{amb1}
\TT(P^n_1, \de_1) \cup \TT(P^n_2, \de_2) \se \CC(E^n) \tx{ for all } n \in \N,
\end{equation}
where $P_1$ and $P_2$ are different probability measures. This situation is excluded by the following proposition:

\bg{prop} If $P_1$ and $P_2$ govern the same experiment $E$, then $P_1=P_2$.
\end{prop}

\bg{proof} Since the classes $\CC(E^n)$ satisfy the same property as the classes $\CC_n$ of Theorem \ref{th}, the proof descends from that theorem.
\end{proof}

The second situation occurs if the same probability measure governs two different experiments. In formal terms, this happens if there is $\de < 1$ such that
\be 
\TT(P^n, \de) \se  \CC(E^n_1) \cap \CC(E^n_2) \tx{ for all } n \in \N,
\ee
where $E_1$ and $E_2$ are two empirically distinguishable experiments. This possibility essentially implies that the classes $\{\CC(E^n)\}_{n \in \N}$ do not determine the whole empirical content of the experiment $E$. Proving or disproving this possibility is arguably a not easy conceptual task, and it is not addressed in this paper.

The third situation occurs if there are experiments for which no probability measure satisfies the conditions \ref{cps} of Cournot's principle. This possibility too is not considered in this paper.


\section{Conclusion \label{conc}}

Cournot's principle relates the mathematics of probability theory to the empirical world of experiments. The original formulation of Cournot's principle applies to classical probability. In this paper, a formulation of the principle suitable for measure-theoretic probability has been proposed. While the classical formulation predicts an empirical consequence from a given probability measure (in this case classical probability), the new formulation is a general criterion for relating probability measures and experiments. Both formulations are based on the correspondence between typical and practically certain events.

In spite of the vagueness of the notions involved, the new version of the principle is formulated in a rigorous manner, and it is rigorously proven that at most one probability measure can be related to an experiment. The possibilities that a probability measure could be related to two different experiments or that no probability measure could be related to an experiment remain open questions.

\typeout{} \bibliography{general}

\end{document}